\documentclass[12pt]{article}
\usepackage[numbers]{natbib}
\usepackage{amssymb, latexsym, amsmath, graphicx, epsfig,
amsfonts, amsthm, bm, float, fullpage, ctable, savefnmark}

\newcommand{\btable}[3]{
\begin{table}[htbp]
\begin{center}
\caption{#2\label{#3}}
\begin{tabular}{#1}}
\newcommand{\etable}{\end{tabular}
\end{center}
\end{table}}

\newtheorem{theorem}{Theorem}
\newtheorem{lemma}{Lemma}
\theoremstyle{definition}
\newtheorem{remark}{Remark}

\newcommand{\mc}[2]{\multicolumn{#1}{c}{#2}}
\newcommand{\bR}{\mathbb{R}}
\newcommand{\bZ}{\mathbb{Z}}
\newcommand{\N}{\mathcal{N}}
\newcommand{\lrp}[1]{\left(#1\right)}
\newcommand{\lrb}[1]{\left\{#1\right\}}
\newcommand{\bt}{\bm{\tau}}

\newenvironment{problist}
   {
      \begin{list}
         {---}
         {
            \setlength{\itemsep}{.5ex}
            \setlength{\parsep}{0ex}
            \setlength{\parskip}{0ex}
            \setlength{\topsep}{.5ex}
         }
   }
   {
      \end{list}
   }

\DeclareMathOperator{\E}{E} \DeclareMathOperator{\var}{var}
 \DeclareMathOperator{\MSE}{MSE}
\DeclareMathOperator{\bias}{bias}

\begin{document}

\bibliographystyle{plainnat}

\date{}

\title{Higher-Order estimation of $s$-th order spectra with
flat-top lag-windows}

\author{Arthur Berg\qquad Dimitris N. Politis\\
{\normalsize {\it University of California, San Diego}} }

\maketitle \abstract{ Improved performance in higher-order spectral
density estimation is achieved using a general class of
infinite-order kernels. These estimates are asymptotically less
biased but with the same order of variance as compared to the
classical estimators with second-order kernels. A simple,
data-dependent algorithm for selecting bandwidth is introduced and
is shown to be consistent with estimating the optimal bandwidth.
Bispectral simulations with several standard models are used to
demonstrate the performance of the proposed methodology.
\vspace{.5cm}

}

\section{Introduction}

Lag-window estimation of the high-order spectra under various
assumptions is known to be consistent and asymptotically normal
\cite{brillinger67a,brillinger67b, lii90cov, lii90norm}.  However,
convergence rates of the estimators depend on the order, or
characteristic exponent, of the lag-window used.  In general,
increasing the order of the lag-window decreases the bias without
affecting the order of magnitude of the variance, thus producing an
estimator with a faster convergence rate. Although estimators using
lag-windows with large orders yield estimates with better mean
square error (MSE) rates, they were overlooked and rarely used in
practice mainly because of two issues. Firstly, in estimating the
second-order spectral density, lag-windows of order larger than two
\emph{may} yield negative estimates, despite the fact that the true
spectral density is known to be nonnegative.  This problem only
pertains, if ever, to the second-order spectral density (since
higher-order spectra are complex-valued), and is easily remedied by
truncating the estimator to zero if it does go negative (thus
improving the already optimal convergence rates \cite{politis95}).
Secondly, when a lag-window has order larger than necessary, the
rate of convergence is still optimal, but the multiplicative
constant will be suboptimal \cite{marron88}.  The second problem is
encountered when using a poor choice of large order lag-window like
the box-shaped truncated lag-window \cite{politis95}, but there are
many other alternatives with descent small-sample performance.
Additionally, when the underlying spectral density is sufficiently
smooth, this second issue is irrelevant since the lag-window with
the largest order performs best. The next section introduces a
family of infinite-order lag-windows for estimating the spectral
density and higher-order spectra.

The use of infinite-order lag-windows is particularly adept to the
estimation of higher-order spectra.  Under the typical scenario of
exponential decay of the autocovariance function (refer to part (ii)
of Theorem 1 within), the MSE rates for estimating the second-order
spectral density using a lag-window of order 2 and an infinite-order
lag-window are  $N^{-4/5}$ and $(\log N)/N$ respectively.  However,
when estimating the third-order spectral density, or bispectrum, the
MSE rates become $N^{-2/3}$ and $(\log N)/N$ respectively.  The
disparity grows stronger with yet higher-order spectra.

The problem of choosing the best bandwidth still remains.   The
optimal bandwidth typically depends on the unknown spectral density
leading to a circular problem--estimation of the spectrum requires
estimation of the bandwidth which in turn requires estimation of the
spectrum. There have been many fixes to this problem; see
\cite{marron96} for a survey of several methods. Section
\ref{section:bandwidth} introduces a new simple, data-dependent
method of determining the bandwidth which is shown to converge to
the asymptotically ideal bandwidth for flat-top lag-windows. An
alternative bandwidth selection algorithm is also included that is
designed for use with second-order lag-windows.  This algorithm uses
the plug-in principle for bandwidth selection but with the flat-top
estimators as the plug-in pilots.

Particular attention is given to the bispectrum as it is a key tool
in several linearity and Gaussianity tests including \cite{Hinich82}
and \cite{gabr80}.  The general bandwidth selection algorithm is
refined and expanded for the bispectrum.  Bispectral simulations
compare two different flat-top lag-windows estimators of the
bispectrum with accompanying bandwidth selection algorithm to the
lag-window estimator using the order two ``optimal'' lag-window and
plug-in bandwidth selection procedure as described in Rao
\cite{rao84}.

We define the flat-top lag-window estimate in Section 2 and derive
its higher-order MSE convergence in Theorem \ref{theorem:bias} under
the ideal bandwidth.  In Section 3, a bandwidth selection algorithm
tailored to the flat-top estimate is introduced and is shown to
automatically adapt to the smoothness of the underlying spectral
density and converge in probability to the ideal bandwidth.  The
focus is then shifted to the bispectrum in Section 4 where the most
general function invariant under the symmetries of the bivariate
cumulant function is constructed.  The bandwidth algorithm is
specialized for the bispectrum, and a separate bandwidth algorithm
for second-order lag-windows is included that is based on the
plug-in method with flat-top estimators as pilots. Simulations of
the bispectrum in Section 5 exhibit the strength of the flat-top
estimators and the bandwidth algorithms.

\section{Asymptotic performance of a general flat-top window}

Let $\bm{x}_1,\bm{x}_2,\ldots, \bm{x}_N$ be a realization of an
$r$-vector valued $s^{\text{th}}$-order stationary (real valued)
time series $\bm{X}_t=(X_t^{(1)},\ldots,X_t^{(r)})'$ with (unknown)
mean $\bm{\mu}=(\mu^{(1)},\ldots,\mu^{(r)})'$. Consider the
$s^{\text{th}}$-order central moment
\begin{equation}
\label{eqn:central}
  C''_{a_1,\ldots,a_s}(\tau_1,\ldots,\tau_{s})=\E\left[(X^{(a_1)}_{t+\tau_1}-\mu^{(a_1)})\cdots
  (X^{(a_s)}_{t+\tau_s}-\mu^{(a_s)})\right].
\end{equation}
where the right-hand side is independent of the choice of $t\in\bZ$.
Stationarity allows us to write the above moment as function of
$s-1$ variables, so we define
$C'_{a_1,\ldots,a_{s}}(\tau_1,\ldots,\tau_{s-1})=C''_{a_1,\ldots,a_{s}}(\tau_1,\ldots,\tau_{s-1},0)$.
For notational convenience, the sequence ${a_1,\ldots,a_s}$ will be
dropped, so $C'_{a_1,\ldots,a_{s}}(\tau_1,\ldots,\tau_{s-1})$ will
be denoted simply by $C'(\bm{\tau})$.  Also $\tau_s$ will
occasionally be used, for convenience, with the understanding that
$\tau_s=0$.

We express the $s^{\text{th}}$-order joint cumulant as
\[
C_{a_1,\ldots,a_{s}}(\tau_1,\ldots,\tau_{s-1})=
\sum_{(\nu_1,\ldots,\nu_p)}(-1)^{p-1}(p-1)!\,\mu_{\nu_1}\cdots\mu_{\nu_p}
\]
where the sum is over all partitions $(\nu_1,\ldots,\nu_p)$ of
$\{0,\ldots,\tau_{s-1}\}$ and $\mu_{\nu_j}=\E\left[\prod_{\tau_i\in
\nu_j}X_{\tau_i}^{(a_{i})}\right]$; refer to \cite{J06} for another
expression of the joint cumulant.  The ($s^\text{th}$-order)
spectral density is defined as
\begin{equation*}
  f(\bm{\omega})=\frac{1}{(2\pi)^{s-1}}
  \sum_{\bm{\tau}\in\bZ^{s-1}}C(\bm{\tau})e^{-i\bt\cdot\bm{\omega}}.
\end{equation*}
We adopt the usual assumption on $C(\bt)$ that it be absolutely
summable, thus guaranteeing the existence and continuity of the
spectral density.  A natural estimator of $C(\bm{\tau})$ is given by
\begin{equation}
\label{eq:Chat} \widehat{C}(\tau_1,\ldots,\tau_{s-1})=
\sum_{(\nu_1,\ldots,\nu_p)}(-1)^{p-1}(p-1)!\,\hat{\mu}_{\nu_1}\cdots\hat{\mu}_{\nu_p}
\end{equation}
where
\[
\hat{\mu}_{\nu_j}=\frac{1}{N-\max(\nu_j)+\min(\nu_j)}\sum_{k=-\min(\nu_j)}^{N-\max(\nu_j)}\,\prod_{t\in\nu_j}x_{t+k}^{(a_j)}
\]

It turns out that the second-order and third-order cumulants, those
that give rise to the spectrum and bispectrum respectively, are
precisely the second-order and third-order central moments
\ref{eqn:central}. Therefore, in these cases, we can greatly
simplify $\hat{C}(\bm{\tau})$ to
\begin{equation}
\label{eq:Chat2}
\widehat{C}(\bm{\tau})=\frac{1}{N}\sum_{t=1}^{N-\gamma}
  \prod_{j=1}^{s}(x^{(a_j)}_{t-\alpha+\tau_j}-\bar{x}^{(a_j)}),
\end{equation}
where $\alpha=\min(0,\tau_1,\ldots,\tau_{k-1})$,
$\gamma=\max(0,\tau_1,\ldots,\tau_{k-1})-\alpha$, and
$\bar{x}^{(a_\ell)}=\frac{1}{N}\sum_{j=1}^Nx^{(a_\ell)}_j$ for
$\ell=1,\ldots s-1$.  We extend the domain of $\hat{C}$ to $\bZ^{s}$
by defining $\hat{C}(\bt)=0$ when then sum in (\ref{eq:Chat}) or
(\ref{eq:Chat2}) is empty.

\noindent Consider a flat-top lag-window function
$\lambda:\bR^{s-1}\rightarrow\bR$ satisfying the following
conditions:
\begin{problist}
  \item[(i)]  $\lambda(\bm{x})\equiv 1$ for all $\bm{x}$ satisfying $\|\bm{x}\|\le
  b$, for some positive number $b$.
  \item[(ii)]  $|\lambda(\bm{x})|\le 1$ for all $\bm{s}$.
  \item[(iii)] For $M\rightarrow\infty$ as $N\rightarrow\infty$,
  but with $M/N\rightarrow 0$,
  \[\lim_{M\rightarrow\infty} \frac{1}{M^{s-1}}\sum_{\|\bm{x}\|\le N}
  \lambda\left(\frac{\bm{x}}{M}\right)<\infty.\]
  \item[(iv)]  $\lambda\in \mathrm{L}_2(\bR^{s-1})$
\end{problist}

 The window $\lambda(x)$ is a ``flat-top'' because of condition
(i); namely, it is constant in a neighborhood of the origin.  The
constant $b$ in (i) is used below in constructing the spectral
density estimate.

Technically, just requiring $\lambda$ just to be bounded could
replace criterion (ii), but there is no benefit in allowing the
window to have values larger than 1. Finally, criteria (iii) and
(iv) are satisfied if, for example, $\lambda$ has compact support.

  Define
$\lambda_M(\bm{t})=\lambda(\bm{t}/M)$ and consider the smoothed
$s^\text{th}$-order periodogram
\begin{equation}
\label{eq:fhat}
  \hat{f}(\bm{\omega})=\frac{1}{(2\pi)^{s-1}}\sum_{\|\bt\|<N}\lambda_M(\bm{\tau})\widehat{C}(\bm{\tau})e^{-i\bt\cdot\bm{\omega}}.
\end{equation}
There is an equivalent expression to this estimator in the frequency
domain given by
\[
\hat{f}(\bm{\omega})=\Lambda_M\ast
I_{a_1,\ldots,a_s}(\bm{\omega})=\int_{\bR^{s-1}}\Lambda_M(\bm{\omega}-\bt)
I_{a_1,\ldots,a_s}(\bt)\,d\bt
\]
where $\Lambda_M$ is the Fourier transform of $\lambda_M$ and
$I_{a_1,\ldots,a_s}$ is the $(s-1)^{\text{th}}$ order periodogram;
namely,
\[
\Lambda_M(\bt)=\int_{\bR^{s-1}} \lambda_M(\bt)e^{-i\bm{\omega}\cdot
\bt}\,d\bt
\]
and
\[
I_{a_1,\ldots,a_s}(\bm{\omega})=\frac{1}{(2\pi)^{s-1}}\sum_{ \bt\in
\bZ^{s-1}}\hat{C}(\bt)e^{-i\bt\cdot\bm{\omega}}
\]
However, equation (\ref{eq:fhat}) is computationally simpler, and it
is this version that will be used throughout the remainder of this
article.

The asymptotic bias convergence rate (and thus the overall MSE
convergence rate) of the estimator (\ref{eq:fhat}) with a flat-top
lag-window $\lambda$ is superior to traditional estimators using
second-order lag-windows. The convergence rates of our estimator
improve with the decay rate of the cumulant function
$C(\bm{\tau})$--the faster the decay to zero, the faster the
convergence.  The following theorem outlines convergence rates under
three scenarios: when the decay of $C(\bt)$ is polynomial,
exponential, and identically zero after some finite time (like an
MA($q$) process). Throughout, conditions on the time series are
assumed so that
\begin{equation}
\label{eq:var}
\var\lrp{\hat{f}(\bm{\omega})}=O\lrp{\frac{M^{s-1}}{N}}.
\end{equation}
This is a very typical assumption and is satisfied under summability
conditions of the cummulants \cite{brillinger67a} or under certain
mixing condition assumptions \cite{lii90norm}.

\begin{theorem}
\label{theorem:bias}  Let $\{\bm{X}_t\}$ be an $r$-vector valued
$s^{\text{th}}$-order stationary time series with unknown mean
$\bm{\mu}$.  Let $\hat{f}(\bm{\omega})$ be the estimator as defined
in (\ref{eq:fhat}) and assume (\ref{eq:var}) is satisfied.

\begin{problist}
\item[(i)]  Assume for some $k\ge 1$,
$\sum_{\bt\in\bZ^{s-1}}\|\bm{\tau}\|^k|C(\bm{\tau})|<\infty$ and
$M\sim aN^c$ with $c=(2k+s-1)^{-1}$, then
\begin{equation}
\label{eq:bias}
\sup_{\bm{\omega}\in[-\pi,\pi]^{s-1}}\left|\bias\lrb{\hat{f}(\bm{\omega})}\right|=o\lrp{N^{\frac{-k}{2k+s-1}}}
\end{equation}
and
\[\MSE(\hat{f}(\bm{\omega}))=O\lrp{N^{\frac{-2k}{2k+s-1}}}.\]
\item[(ii)]  Assume $C(\bm{\tau})$ decreases
  geometrically fast, i.e. $|C(\bm{\tau})|\le De^{-d\|\bm{\tau}\|}$, for
  some positive constants $d$ and $D$ and $M\sim A\log N$ where $A\ge1/(2db)$, then
\begin{equation}
 \label{eq:bias2}
\sup_{\bm{\omega}\in[-\pi,\pi]^{s-1}}\left|\bias\lrb{\hat{f}(\bm{\omega})}\right|=
O\lrp{\frac{1}{\sqrt{N}}}
\end{equation}
and
  \begin{equation}
  \label{eq:mse2}
  \MSE\lrp{\hat{f}(\bm{\omega})}=O\lrp{\frac{\log N}{N}}.
  \end{equation}

\item[(iii)]  Assume $C(\bm{\tau})=0$ for $\|\bt\|>q$ and let $M$ be constant such that $bM\ge q$,
then
\[
\sup_{\bm{\omega}\in[-\pi,\pi]^{s-1}}\left|\bias\lrb{\hat{f}(\bm{\omega})}\right|=O\lrp{\frac{1}{N}}
\]
and
\[
\MSE\lrp{\hat{f}(\bm{\omega})}=O\lrp{\frac{1}{N}}
\]
\end{problist}
\end{theorem}
\begin{remark}
Equations (\ref{eq:bias}), respectively (\ref{eq:bias2}), remain
true with the assumptions on $M$ replaced with
$M^{k+s-1}/N\rightarrow 0$, respectively $e^MM^{s-1}/N\rightarrow
  0$.
\end{remark}

\begin{remark}
Depending on the constant $A$ in part (ii), the bias in
(\ref{eq:bias2}) may be as small as $O\lrp{(\log N)^{s-1}/N}$.
\end{remark}

\begin{remark}
We do not assume the mean $\bm{\mu}$ of the time series is known.
This adds an extra term of order $O(M^{s-1}/N)$ to the bias; see the
proof of Theorem \ref{theorem:bias} in the appendix for further
details.
\end{remark}

\begin{remark}
Traditional estimators using second-order lag-windows have bias
convergence rates of order $O(1/M^2)$ regardless of the three
scenarios listed in Theorem \ref{theorem:bias}.  However when the
spectral density is smooth enough, like in the case of an ARMA
process (where $C(\bt)$ decays exponentially), traditional
estimators perform considerably worse.  For example, estimation of
the bispectrum of an ARMA process has an asymptotic MSE rate of
$N^{-2/3}$ in the traditional case, but an asymptotic MSE of $(\log
N)/N$ using flat-top lag-windows. The distinction is even more
profound in estimating higher-order spectra where the best rate
achieved is $N^{-4/(3+s)}$ for traditional estimators and again
$(\log N)/N$ using flat-top lag-windows.  Even in the worst case of
polynomial decay, our proposed estimator still beats, or possibly
ties with, traditional estimators in terms of asymptotic MSE rates.
\end{remark}

The asymptotic analysis in Theorem \ref{theorem:bias} relies on
having the appropriate bandwidth $M$ based on the various decay
rates $C(\bt)$. In the next section we propose an algorithm that,
for the most part, automatically detects the correct decay rate of
$C(\bt)$ and supplies the practitioner with an asymptotically
consistent estimate of $M$.

\section{A Bandwidth Selection Procedure}
\label{section:bandwidth}

For $\bm{\tau}\in\bZ^{s-1}$, consider the normalized cumulant
function
\[\rho(\bt)=\frac{C(\bt)}{\left(\prod_{i=1}^sC_{a_i}(0)\right)^{1/2}}\]
with natural estimator
\[\hat{\rho}(\bt)=
\frac{\hat{C}(\bt)}{\left(\prod_{i=1}^s\hat{C}_{a_i}(0)\right)^{1/2}}.\]
Let $B_{x,y}$ ($x,y>0$) denote the set of indices in $\bZ^{s-1}$
contained in the half-open $s-1$-dimensional annulus of inner radius
$x$ and outer radius $y$, i.e.
\begin{equation}
\label{eq:norm} B_{x,y}=\{\bt\in\bZ^{s-1}:x<\|\bt\|\le y\}.
\end{equation}

The following algorithm for estimating the bandwidth of a flat-top
estimator is a multivariate extension of an algorithm proposed in
\cite{politis03}.

\begin{quote}
  \textsc{Bandwidth Selection Algorithm}\\
  Let $k>0$ be a fixed constant, and $a_N$ be a nondecreasing
  sequence of positive integers tending to infinity such that $a_N=o(\log
  N)$.
  Let $\hat{m}$ be the smallest number such that
  \begin{equation}
  \label{eq:alg}
|\hat{\rho}(\bt)|<k\sqrt{\frac{\log_{10} N}{N}}\qquad\text{ for all
}\bt\in B_{\hat{m},\hat{m}+a_N}
  \end{equation}
Then let $\hat{M}=\hat{m}/b$ (where $b$ is the ``flat-top radius''
as defined by condition (i) of a flat-top lag-window).
\end{quote}
\begin{remark}
A norm was not specified in (\ref{eq:norm}) and any norm may be
used.  The sup norm, for example, may be preferable to the Euclidean
norm in practice since the region in (\ref{eq:norm}) becomes
rectangular instead of circular.
\end{remark}

\begin{remark}
\label{rmk:subsample} The positive constant $k$ is irrelevant in the
asymptotic theory, but is relevant for finite-sample calculations.
In order to determine an appropriate value of $c$ for computation,
we consider the following approximation
\begin{equation}
\label{eq:approxNorm}
\sqrt{N}\left(\hat{\rho}(\bt_0)-\rho(\bt_0)\right)\, \dot{\sim}\,
\N\left(0,\sigma^2\right).
\end{equation}
This approximation holds under general assumptions of the time
series and for any fixed  $\bm{\tau}_0\in\bZ^{s-1}$.  The variance
$\sigma^2$ does not depend on the choice of $\bm{\tau}_0$ provided
$\bm{\tau}_0$ is not a ``boundary point''; see \cite{brillinger67b}
for more details.  Let $\hat{\sigma}$ be the estimate of $\sigma$
via a resampling scheme like the block bootstrap.  A approximate
pointwise 95\% confidence bound for $\rho(\cdot)$ is given by
$\frac{\pm 1.96\, \hat{\sigma}}{\sqrt{N}}$. Therefore if we let
$a_N=5$, then $k=2\hat{\sigma}$ generates an approximate 95\%
simultaneous confidence bound by Bonferroni's inequality by noting
that $\sqrt{\log_{10} N}\approx 1.5$ for moderately sized $N$.

\end{remark}

The bandwidth selected using the above procedure converges precisely
to the ideal bandwidth in each of the three cases of Theorem
\ref{theorem:bias}, as is proved in the following theorem under the
two natural assumptions in (\ref{eq:rho1}) and (\ref{eq:rho2})
below\footnote{Under general regularity conditions, (\ref{eq:rho1})
holds as does the even stronger assumption of $\sqrt{N}$ asymptotic
normality, and (\ref{eq:rho2}) holds from general theory of extremes
of dependent sequences; refer to \cite{Leadbetter83}}.

\begin{theorem}
\label{theorem:bandwidth} Assume conditions strong enough to ensure
that for any fixed $n$,
\begin{equation}
\label{eq:rho1} \max_{\bt\in
B_{0,n}}|\hat{\rho}(\bm{\sigma}+\bt)-\rho(\bm{\sigma}+\bt)|=
O_p\left(\frac{1}{\sqrt{N}}\right)
\end{equation}
uniformly in $\bm{\sigma}$, and for any $M$, that may depend on $N$,
the following holds
\begin{equation}
\label{eq:rho2} \max_{\bt\in
B_{0,M}}|\hat{\rho}(\bm{\sigma}+\bt)-\rho(\bm{\sigma}+\bt)|=
O_p\left(\sqrt{\frac{\log M}{N}}\right)
\end{equation}
uniformly in $\bm{\sigma}$.
\begin{problist}
  \item[(i)]  Assume $C(\bt)\sim A\|\bt\|^{-d}$
  for some positive constants $A$ and $d\ge1$. Then
  \[
\hat{M}\stackrel{P}{\sim}A_0\frac{N^{1/2d}}{(\log N)^{1/2d}}
  \]
where $A_0=A^{1/d}/(k^{1/d}b)$; here $A\stackrel{P}{\sim}B$ means
$A/B\rightarrow 1$ in probability.
\item[(ii)]
Assume $C(\bt)\sim A\,\xi^{\|\bt\|}$ for some positive constant $A$
and $|\xi|<1$. Then
  \[
\hat{M}\stackrel{P}{\sim}A_1\log N
  \]
where $A_1=-1/(b\log|\xi|)$.
\item[(iii)]  Suppose $C(\bt)=0$ when $\|\bt\|>q$, but $C(\bt)\not=0$
for some $\bt$ with norm $q$,~then~$\hat{M}~\stackrel{P}{\sim}~q/b$.
\end{problist}
\end{theorem}

\section{Bispectrum}
\label{bisp}

Now we will focus on estimating the bispectrum using flat-top
lag-windows.  The third-order cumulant reduces to the third-order
central moment with estimator given by (\ref{eq:Chat2}). It is
easily seen that the third-order central moment, $C(\tau_1,\tau_2)$,
satisfies the following symmetry relations:
\begin{equation}
\label{eq:symmetry}
  C(\tau_1,\tau_2)=C(\tau_2,\tau_1)=C(-\tau_1,\tau_2-\tau_1)=
  C(\tau_1-\tau_2,-\tau_2)
\end{equation}
Naturally, we would expect the lag-window function,
$\lambda(\tau_1,\tau_2)$, in the estimator (\ref{eq:fhat}), to
posses the same symmetries.  So if a lag-window $\lambda$ does not a
priori have the symmetries as in (\ref{eq:symmetry}), we can
construct a symmetrized version given by
\begin{equation}
\label{eq:symmetrize}
\tilde{\lambda}=g\left(\lambda(x,y),\lambda(y,x),\lambda(-x,y-x),
\lambda(y-x,-x),\lambda(x-y,-y),\lambda(-y,x-y)\right)
\end{equation}
where $g$ is any symmetric function (of its six variables); for
example $g$ could be the geometric or arithmetic mean.  It is worth
noting that the symmetrized version of $\lambda$ is connected to the
theory of group representations of the symmetric group $S_3$.  As a
special case, symmetric lag-windows can be constructed from a
one-dimensional lag-window $\lambda(x)$, namely,
\begin{equation}
\label{eq:symmetrize2}
\tilde{\lambda}=g\left(\lambda(x),\lambda(y),\lambda(-x),
\lambda(y-x),\lambda(x-y),\lambda(-y)\right)
\end{equation}
and if $\lambda(x)$ is an even function, then (\ref{eq:symmetrize2})
becomes $\tilde{\lambda}=h\left(\lambda(x),\lambda(y),
\lambda(y-x)\right)$ where $h$ is any symmetric function (of its
three variables).

 Several choices of lag-windows are considered in
\cite{rao03} including the so-called ``optimal window'',
$\lambda_{\text{opt}}$, which is in some sense optimal among
lag-windows of order 2; see Theorem 2 on page 43 of \cite{rao84}.
This lag-window is defined as \cite{Saito85}
\[
\lambda_{\text{opt}}(\tau_1,\tau_2)=\frac{8}{\alpha(\tau_1,\tau_2)^2}
J_2(\alpha(\tau_1,\tau_2))
\]
where $J_2$ is the second-order Bessel function of the first kind,
and
\[
\alpha(x,y)=\frac{2\pi}{\sqrt{3}}\sqrt{x^2-xy+y^2}
\]
Although $\lambda_{\text{opt}}$ is optimal among order 2
lag-windows, it is sub-optimal to higher-order lag-windows, such as
flat-top lag-windows. Also, since $\lambda_{\text{opt}}$ is not
compactly supported, it has the potential of being computationally
taxing.

We detail two simple flat-top lag-windows satisfying the symmetries
in (\ref{eq:symmetry}), but the supply of examples is limitless by
(\ref{eq:symmetrize}). The first example is a right pyramidal
frustum with the hexagonal base $|x|+|y|+|x-y|=2$. We let
$c\in(0,1)$ be the scaling parameter that dictates when the frustum
becomes flat, that is, the flat-top boundary is given by
$|x|+|y|+|x-y|=2c$. The equation of this lag-window is given by
\[
\lambda_{\text{rpf}}(\tau_1,\tau_2)=\frac{1}{1-c}
\lambda_{\text{rp}}(\tau_1,\tau_2)-\frac{c}{1-c}
\lambda_{\text{rp}}\left(\frac{\tau_1}{c},\frac{\tau_2}{c}\right)
\]
where $\lambda_{\text{p}}$ is the equation of the right pyramid with
base $|x|+|y|+|x-y|=2$, i.e.,
\[
\lambda_{\text{rp}}(x,y)=
\begin{cases}
(1-\max(|x|,|y|))^+  ,&-1\le x,y\le 0\ \text{or}\ 0\le x,y\le 1\\
(1-\max(|x+y|,|x-y|))^+, &\text{otherwise}
\end{cases}
\]

 The second flat-top lag-window that we propose is the right
conical frustum with elliptical base $x^2-xy+y^2=1$. As in the
previous example, there is a scaling parameter $c\in(0,1)$, and the
lag-window becomes flat in the ellipse $x^2-xy+y^2=c^2$.  The
equation of this lag-window is given by
\[
\lambda_{\text{rcf}}(\tau_1,\tau_2)=\frac{1}{1-c}
\lambda_{\text{rc}}(\tau_1,\tau_2)-\frac{c}{1-c}
\lambda_{\text{rc}}\left(\frac{\tau_1}{c},\frac{\tau_2}{c}\right)
\]
where $\lambda_{\text{rc}}$ is the equation of the right cone with
base $x^2-xy+y^2=1$, i.e.,
\[
\lambda_{\text{rc}}(x,y)= (1-\sqrt{x^2-xy+y^2})^+
\]

Although in both examples the value for $b$, as defined in property
(i) of the flat-top lag-window function, is smaller than the
parameter $c$, the symmetries (\ref{eq:symmetry}) permit us to only
consider the region $0\le y\le x$ for which a circular arc of radius
$c$ does fit. So in the two examples above, we take the value of $b$
to be the parameter $c$.

\begin{figure}
\centerline{\includegraphics[width=6in]{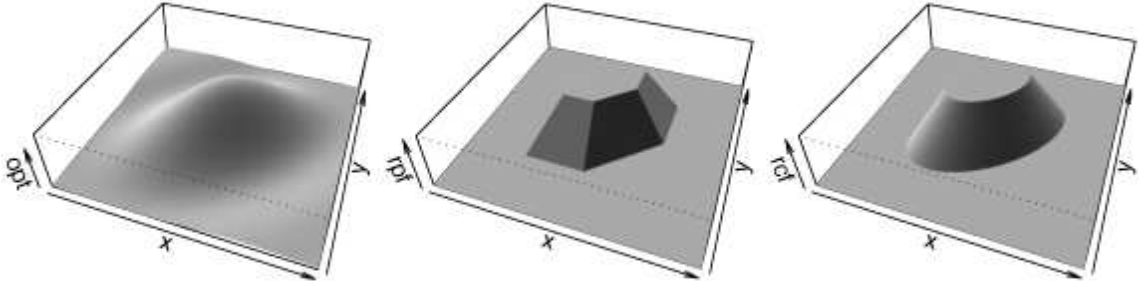}} \caption{Plots
of the three lag-windows, $\lambda_{\text{opt}}$,
$\lambda_{\text{rpf}}$, and $\lambda_{\text{rpc}}$ (with $c=1/2$ in
the latter two).}
\end{figure}

The bandwidth selection algorithm can be refined in the context of
the bispectrum.  The symmetries in (\ref{eq:symmetry}) allow
restriction to the region
\begin{equation}
\label{sector1} \{(\tau_1,\tau_2)\in\bR^2\ |\ 0\le \tau_2\le\tau_1\}
\end{equation}
Here is the modified bandwidth selection algorithm for flat-top
kernels that is tailored to the bispectrum:
\begin{quote}
  \textsc{Practical Bandwidth Selection Algorithm for the Bispectrum}\\
Let $\tilde{k}=k_1>0$ if $n=1$, otherwise $\tilde{k}=k_2>0$, and let
$L$ be a positive integer that is $o\left(\log N\right)$.  Order the
points $\{(\tau_1, \tau_2)\in\bZ^2\ |\ 0 <
\tau_2<\tau_1\}\cup\{(1,0)\}$ with the usual lexicographical
ordering, so $P_1=(1,0)$, $P_2=(2,1)$, $P_3=(3,1)$, $P_4=(3,2)$, and
so forth; in general, $P_n=(i,j)$ where
$i=\left\lfloor\left(\frac{3}{2}+\sqrt{2n-2}\right)\right\rfloor$
and $j=n-\frac{1}{2}\left(i^2-3i\right)-2$.  Let $\hat{m}$ be the
smallest number such that
\begin{equation}
\left|\hat{\rho}\left(P_{\hat{m}+\ell}\right)\right|<\tilde{k}\sqrt{\frac{\log
N}{N}}\qquad\text{ for all }\ell=1,\ldots,L.
  \end{equation}
Then let $\hat{M}=\left(\text{first coordinate of
}P_{\hat{m}}\right)/b=
\left(\left\lfloor\left(\frac{3}{2}+\sqrt{2\hat{m}-2}\right)\right\rfloor\right)/b$.
\end{quote}

\begin{remark}
Except for the first point, $(1,0)$, this algorithm does not
incorporate boundary points since the asymptotic variance is larger
on the boundary;  the first point is included as there are no
interior points with first coordinate equal to 1.  The constant
$\tilde{k}$ is adjusted to account for the larger variance in the
first point by providing a separate threshold, $k_1$, for this
point.
\end{remark}
\begin{remark}
\label{rmk:subsample2} As suggested with the general algorithm, a
subsampling procedure should be used to determine the appropriate
constants $k_1$ and $k_2$.  However, one should be careful when
choosing a point $\bm{\tau}_0$ for the approximation
(\ref{eq:approxNorm}) since high variances at the origin and on the
boundary tend to cause high variances near the origin and near the
boundary in finite-sample scenarios.  Therefore an interior point
like $(6,3)$ (as opposed to $(2,1)$) should be used in determining
$k_2$, and a point like (3,0) (as opposed to $(1,0)$) should be used
in determining $k_1$.
\end{remark}

A modified bandwidth selection procedure is now proposed for use
with the sub-optimal lag-windows of order 2. In this case, we
propose using a bandwidth selection procedure based on the usual
``solve-the-equation plug-in" approach \cite{marron96}, but with
flat-top estimates of the unknown quantities as the plug-in pilots.
This will afford faster convergence rates of the bandwidth as
compared estimates based on second-order pilots as well as solve the
problem of selecting bandwidths for the pilots.

The optimal bandwidth at each point in the region (\ref{sector1}),
when using differentiable second order kernels, is derived in
\cite{rao84}, and is given by
\begin{equation}
\label{eq:optimalM}
\begin{split}
 M_\lambda(\omega_1,\omega_2)=&\Bigg\{\frac{\pi N\,
}{\|\lambda\|_{L_2} f(\omega_1)f(\omega_2)f(\omega_1+\omega_2)}
\left(\frac{\partial^2\lambda(\tau_1,\tau_2)}{\partial
\tau_1\,\partial\tau_1}\bigg|_{\tau_1=\tau_2=0}\right)^2\times\\
&\quad\times\left|\left(\frac{\partial^2}{\partial
\omega_1^2}-\frac{\partial^2}{\partial\omega_1\partial\omega_2}
+\frac{\partial^2}{\partial \omega_2^2}\right)
f(\omega_1,\omega_2)\right|^{2}\Bigg\}^{\frac{1}{6}}
\end{split}
\end{equation}

Estimates of the spectral density using flat-top lag-windows is
discussed above, and estimating the partial derivatives of the
bispectrum follow similarly.  For instance, the three second order
partial derivatives needed in (\ref{eq:optimalM}) can be estimated
by

\begin{equation}
\label{eq:partialfhat}
\begin{split}
  \widehat{f_{\omega_i,\omega_j}}(\omega_1,\omega_2)&=
  \frac{\partial^2}{\partial
  \omega_i\partial\omega_j}\hat{f}(\omega_1,\omega_2)\\
  &=\frac{1}{(2\pi)^{2}}
  \sum_{\tau_1=-N}^N\sum_{\tau_2=-N}^N\tau_i\tau_j\,
  \lambda_M(\tau_1,\tau_2)\widehat{C}(\tau_1,\tau_2)
  e^{-i\bt\cdot\bm{\omega}}\qquad i,j=1,2.
\end{split}
\end{equation}

By mimicking the proof of Theorem \ref{theorem:bias}, the estimator
in (\ref{eq:partialfhat}) has the same asymptotic performance as the
estimator $\hat{f}(\bm{\omega})$ in Theorem \ref{theorem:bias} but
under a slightly stronger assumption for part (i) that
$\sum_{\bt\in\bZ^{2}}\|\bm{\tau}\|^{k+2}|C(\bm{\tau})|<\infty$. We
construct the estimator $\hat{M}_\lambda$ by replacing the unknown
$f$ and its derivatives in (\ref{eq:optimalM}) with flat-top
estimates producing
\[
\begin{split}
 \widehat{M_\lambda}(\omega_1,\omega_2)=&\Bigg\{\frac{\pi N\,
}{\|\lambda\|_{L_2}
\hat{f}(\omega_1)\hat{f}(\omega_2)\hat{f}(\omega_1+\omega_2)}
\left(\frac{\partial^2\lambda(\tau_1,\tau_2)}{\partial
\tau_1\,\partial\tau_1}\bigg|_{\tau_1=\tau_2=0}\right)^2\times\\
&\quad\times\left|\left(\frac{\partial^2}{\partial
\omega_1^2}-\frac{\partial^2}{\partial\omega_1\partial\omega_2}
+\frac{\partial^2}{\partial \omega_2^2}\right)
\hat{f}(\omega_1,\omega_2)\right|^{2}\Bigg\}^{\frac{1}{6}}
\end{split}
\]

The next theorem provides convergence rates of the plug-in algorithm
with flat-top pilots.

\begin{theorem}
\label{theorem:pilot} Assume conditions on $\hat{\rho}$ such that
(\ref{eq:rho1}) and (\ref{eq:rho2}) of Theorem
\ref{theorem:bandwidth} hold true, and assume conditions strong
enough to ensure\footnote{Certain mixing condition assumptions
guarantee this; see \cite{politis03} for an example.}
\[
\var\left(\widehat{f_{\omega_i,\omega_j}}\right)=O\left(\frac{M^{s-1}}{N}\right)\qquad
(i,j=1,2)
\]

\begin{problist}
  \item[(i)]  Assume $C(\bt)\sim A\|\bt\|^{-d}$
  for some positive constants $A$ and $d>s+2$. Then
  \[
\widehat{M_\lambda}=M_\lambda\left(1+O_p\left(\left(\frac{\log
N}{N}\right)^{\frac{\lceil d-s-2\rceil}{2d}}\right)\right).
  \]
\item[(ii)]
Assume $C(\bt)\sim A\,\xi^{\|\bt\|}$ for some positive constant $A$
and $|\xi|<1$. Then
  \[
\widehat{M_\lambda}=M_\lambda\left(1 +O_p\left(\left(\frac{\log
N}{N}\right)^{\frac{1}{2}}\right)\right).
  \]
\item[(iii)]  Suppose $C(\bt)=0$ when $\|\bt\|>q$, but $C(\bt)\not=0$
for some $\bt$ with norm $q$, then
  \[
\widehat{M_\lambda}=M_\lambda\left(1 +O_p\left(
\frac{1}{\sqrt{N}}\right)\right).
  \]
\end{problist}
\end{theorem}

In many cases, the convergence is a significant improvement over the
traditional plug-in approach with second-order lag-window pilots.
For example, the convergence of the bandwidth for data from an ARMA
process would be $M(1+O_P(N^{-2/9}))$ using second-order pilots and
techniques similar to \cite{B93,B96}, but by using flat-top pilots,
the convergence improves to $M(1+O_P(\sqrt{\log N/ N}))$.

\section{Bispectral Simulations}

The three lag-windows detailed above--$\lambda_{\text{opt}}$,
$\lambda_{\text{rpf}}$, and $\lambda_{\text{rcf}}$--are compared by
their mean square error performance in estimating the bispectrum of
four standard time series models.   Three criteria are used to
evaluate the performance of the bispectral estimates. The first two
criteria are the estimators performance in estimating the bispectrum
at the two points $(0,0)$ and $(2,1)$.  The  bispectrum at the point
$(0,0)$ is real-valued, and estimates typically have variances
significantly larger than estimates at the interior point (2,1)
(exactly 30-times larger, asymptotically, if the second-order
spectrum is flat). The bispectrum at the point $(2,1)$ is complex
valued and performance is evaluated based on the estimation of the
real part, complex part, and absolute value.  The third criteria of
evaluation is a composite evaluation of performance of the
estimators over a rough grid of six points, standardized
appropriately (further details below). The simulations are computed
with data from the four stationary time series models: iid
$\chi_1^2$, ARMA(1,1), GARCH(1,1), and bilinear(1,0,1,1).  The first
two are linear time series models whereas the last two nonlinear
models. Two sample sizes, $N=200$ and $N=2000$, are used throughout.
Every simulation is repeated over 500 realizations.

The third criteria of evaluation, the composite evaluation is now
described in further detail.  The symmetries of $C$ as given in
(\ref{eq:symmetry}) induce the following symmetries in the spectral
density:
\[
f(\omega_1,\omega_2)=f(\omega_2,\omega_1)=f(\omega_1,-\omega_1-\omega_2)=
f(-\omega_1-\omega_2,\omega_2)=f^*(-\omega_1,-\omega_2)
\]
The above symmetries in combination with the periodicity of $f$
imply that $f$ can be determined over the entire plane just by its
values in the closed triangle $T$ with vertices $(0,0)$, $(\pi,0)$,
and $(2\pi/3,2\pi/3)$.  So $f$ is estimated at ${n-1\choose
2}=\frac{(n-1)(n-2)}{2}$ equally spaced points inside $T$ with
coordinates $\bm{\omega}_{ij}=\left(\frac{\pi(2i+2j)}{3
n},\frac{2\pi j}{3n}\right)$ where $i=1,\ldots,n-1$ and
$j=1,\ldots,n-i-1$ (we take $n=5$ in the simulations).

The estimates at $\bm{\omega}_{ij}$ are standardized to make them
comparable. Since, for $(\omega_1,\omega_2)$ inside $T$,
\cite{rao84}
\[
\var\left(\hat{f}(\omega_1,\omega_2)\right)\approx \frac{M^2}{
N}\frac{\|\lambda\|_{L_2}}{2\pi}f(\omega_1)f(\omega_2)f(\omega_1+\omega_2),
\]
$\hat{f}(\omega_1,\omega_2)$ is standardized by dividing it by
$\sqrt{f(\omega_1)f(\omega_2)f(\omega_1+\omega_2)}$.  This leads to
the composite evaluation of $\hat{f}$ over a course grid of points
by the quantity
\[
\texttt{err}(\lambda)\triangleq\sum_{i=1}^{n-1}\ \sum_{j=1}^{n-i-1}
\left|\frac{\hat{f}(\bm{\omega}_{ij})-f(\bm{\omega}_{ij})}
{\sqrt{f(\omega^{(1)}_{ij})
f(\omega^{(2)}_{ij})f(\omega^{(1)}_{ij}+\omega^{(2)}_{ij})}}\right|
\]
and the empirical MSE is calculated by averaging
$\texttt{err}(\lambda)^2$ over the 500 realizations.

In the tables of MSE estimates below, the first two rows are
estimates from the flat-top lag-windows $\lambda_{\text{rpf}}$ and
$\lambda_{\text{rcf}}$ with the bandwidth derived from the Bandwidth
Selection Algorithm for the Bispectrum, as described above, with
parameters $L=5$, $c=.51$, and $k$ determined via the block
bootstrap (see Remarks \ref{rmk:subsample} and
\ref{rmk:subsample2}).  The third and fourth rows are estimates
using the $\lambda_{\text{opt}}$ with bandwidths from the plug-in
method with flat-top pilots (f.p.) and second-order pilots (s.p.)
respectively. The first column of each table concerns the estimation
of the bispectrum at $(0,0)$, taking absolute values if the estimate
is complex valued. The next three columns concern the estimation of
the real part, complex part, and absolute value of the bispectrum,
respectively, at the point $(2,1)$.  The last column, labeled $T_6$,
concerns the composite evaluation over a coarse grid of 6 points.

Simulations (based on 1000 realizations) were conducted to determine
the optimal finite-sample bandwidth with minimal MSE (checking up to
a bandwidth size of 20).  In the first three models--IID, ARMA, and
GARCH--the optimal bandwidth is 1 under each evaluation criterion
and every lag-window.  The estimators with best MSE performance in
these models were the estimators with the best bandwidth selection
procedure (the choice of lag-window was somewhat secondary). The
bilinear model, however, had different optimal bandwidths depending
on the evaluation criterion and the lag-window. The optimal
bandwidths for the bilinear model were incorporated into MSE tables
by subscripting each value with the best bandwidth followed by the
second best bandwidth.  The optimality of the flat-top lag-window,
independent of the bandwidth selection procedure, can be observed in
this model as the optimal bandwidths are larger than 1.

Simulations are also carried out to study the bandwidth selection
procedure for the bispectrum.  Histograms, placed in Appendix
\ref{app:hist}, depict the selected bandwidths for each model over
500 realizations under five procedures (a)--(e) described below.
Procedure (a) produces bandwidths for flat-top lag-windows
$\lambda_{\text{rpf}}$ and $\lambda_{\text{rcf}}$ whereas procedures
(b) through (e) produce bandwidths for $\lambda_{\text{opt}}$.
\begin{enumerate}
\item[(a)]  Practical bandwidth selection algorithm for the bispectrum of Section
\ref{bisp}
\item[(b)]  Plug-in method at the origin with flat-top pilots
\footnote{The pilot estimates were derived from the flat-top
lag-windows $\lambda_{\text{rpf}}$ and the trapezoidal flat-top
window \cite{politis95}.  The bandwidths for the pilot estimators
are derived from the bandwidth selection algorithm of Section
\ref{section:bandwidth}.}\saveFN\fna
\item[(c)]  Plug-in method at the point (2,1) with flat-top
pilots\useFN\fna
\item[(d)]  Plug-in method at the origin with second-order pilots \footnote{The
Parzen and optimal lag-windows were used as pilots with bandwidths
$\lfloor N^{1/5}\rfloor$ and $\lfloor N^{1/6}\rfloor$
respectively.}\saveFN\fnb
\item[(e)]  Plug-in method at the point (2,1) with second-order
pilots\useFN\fnb
\end{enumerate}

The performance of the above bandwidth selections procedures are
evaluated by computing MSE estimates based on the simulations
determining the optimal bandwidth.  Since procedure (a) produces a
global bandwidth, comparison is not so straightforward in the
bilinear case where the optimal bandwidth at the origin is different
from that of the interior.

\subsection{IID Data}

Identical and independent $\chi_1^2$ data is generated with a
central third moment $\mu_3=8$.  Therefore the true bispectrum is
$f(\omega_1,\omega_2)\equiv \frac{\mu_3}{(2\pi)^2}\approx .202642$.
  The following tables give the empirical MSE calculations of the
estimated bispectrum over lengths $N=200$ and $N=2000$ based on 500
simulations.
\vspace{.3cm}
\begin{table}[H]
\centering
\begin{tabular}{|c|ccccc|}
  \hline
  $N=200$ & $|\hat{f}(0,0)|$ & Re$\hat{f}(2,1)$ & Im$\hat{f}(2,1)$ & $|\hat{f}(2,1)|$ & $T_6$ \\
  \hline
$\lambda_{\text{rpf}}$ & $\text{0.02796}$ & $\text{0.02061}$ & $\text{3.131e-04}$ & $\text{0.02093}$ & $\text{709.4}$ \\
$\lambda_{\text{rcf}}$ & $\text{0.02778}$ & $\text{0.02060}$ & $\text{3.314e-04}$ & $\text{0.02094}$ & $\text{709.4}$ \\
$\lambda_{\text{opt}} $ (f.p.) & $\text{0.02582}$ & $\text{0.02086}$ & $\text{3.577e-04}$ & $\text{0.02122}$ & $\text{709.8}$ \\
$\lambda_{\text{opt}} $ (s.p.) & $\text{0.02806}$ & $\text{0.02116}$ & $\text{7.121e-04}$ & $\text{0.02187}$ & $\text{715.5}$ \\
  \hline
  \hline
  $N=2000$ & $|\hat{f}(0,0)|$ & Re$\hat{f}(2,1)$ & Im$\hat{f}(2,1)$ & $|\hat{f}(2,1)|$ & $T_6$ \\
  \hline
$\lambda_{\text{rpf}}$ & $\text{2.887e-03}$ & $\text{2.063e-03}$ & $\text{1.799e-05}$ & $\text{2.081e-03}$ & $\text{71.19}$ \\
$\lambda_{\text{rcf}}$ & $\text{2.865e-03}$ & $\text{2.064e-03}$ & $\text{1.875e-05}$ & $\text{2.083e-03}$ & $\text{71.22}$ \\
$\lambda_{\text{opt}} $ (f.p.) & $\text{2.616e-03}$ & $\text{2.101e-03}$ & $\text{2.085e-05}$ & $\text{2.121e-03}$ & $\text{71.23}$ \\
$\lambda_{\text{opt}} $ (s.p.) & $\text{3.294e-03}$ & $\text{2.184e-03}$ & $\text{1.039e-04}$ & $\text{2.288e-03}$ & $\text{71.45}$ \\
  \hline
\end{tabular}
\caption{MSE estimates based on iid data for $N=200$ and $N=2000$.}
\label{tab:iid1}
\end{table}

The flat-top estimators and $\lambda_{\text{opt}}$ (f.p.) outperform
$\lambda_{\text{opt}}$ (f.p.) in every criterion considered. For
$N=2000$, bandwidth procedures (a), (b), and (c) perform extremely
well (refer to the histograms in Figure \ref{fig:iid} in Appendix
\ref{app:hist}) producing the optimal bandwidth 1 over 95\% of the
time in each case.

\subsection{ARMA Model}

The ARMA(1,1) model
\[
X_t=.5X_{t-1}-.5Z_{t-1}+Z_t
\]
is now considered where $Z_t\stackrel{\text{iid}}{\sim}\N(0,1)$.
This time series is Gaussian, so both the bispectrum and normalized
bispectrum are identically zero.

\begin{table}[H]
\centering
\begin{tabular}{|c|lllll|}
  \hline
  $N=200$ & $|\hat{f}(0,0)|$ & Re$\hat{f}(2,1)$ & Im$\hat{f}(2,1)$ & $|\hat{f}(2,1)|$ & $T_6$ \\
  \hline
$\lambda_{\text{rpf}}$ & $\text{6.102e-05}$ & $\text{2.329e-05}$ & $\text{4.468e-06}$ & $\text{2.776e-05}$ & $\text{313.3}$ \\
$\lambda_{\text{rcf}}$ & $\text{6.760e-05}$ & $\text{2.435e-05}$ & $\text{4.624e-06}$ & $\text{2.897e-05}$ & $\text{316.5}$ \\
$\lambda_{\text{opt}} $ (f.p.) & $\text{4.422e-05}$ & $\text{2.172e-05}$ & $\text{5.235e-06}$ & $\text{2.696e-05}$ & $\text{302.8}$ \\
$\lambda_{\text{opt}} $ (s.p.) & $\text{1.198e-04}$ & $\text{3.088e-05}$ & $\text{2.982e-05}$ & $\text{6.070e-05}$ & $\text{412.0}$ \\
  \hline
  \hline
  $N=2000$ & $|\hat{f}(0,0)|$ & Re$\hat{f}(2,1)$ & Im$\hat{f}(2,1)$ & $|\hat{f}(2,1)|$ & $T_6$ \\
  \hline
$\lambda_{\text{rpf}}$ & $\text{2.997e-06}$ & $\text{2.096e-06}$ & $\text{6.896e-08}$ & $\text{2.165e-06}$ & $\text{24.21}$ \\
$\lambda_{\text{rcf}}$ & $\text{3.297e-06}$ & $\text{2.137e-06}$ & $\text{7.359e-08}$ & $\text{2.210e-06}$ & $\text{24.59}$ \\
$\lambda_{\text{opt}} $ (f.p.) & $\text{3.129e-06}$ & $\text{2.132e-06}$ & $\text{2.796e-07}$ & $\text{2.412e-06}$ & $\text{24.74}$ \\
$\lambda_{\text{opt}} $ (s.p.) & $\text{2.142e-05}$ & $\text{4.222e-06}$ & $\text{4.349e-06}$ & $\text{8.571e-06}$ & $\text{33.53}$ \\
  \hline
\end{tabular}
\caption{MSE estimates based on arma data for $N=200$ and $N=2000$.}
\end{table}

The flat-top estimators and $\lambda_{\text{opt}}$ (f.p.) even more
significantly outperform $\lambda_{\text{opt}}$ (f.p.) in this model
for every criterion considered.  Good performance is mostly
attributed to good bandwidth selection, but true optimal properties
of the flat-top lag-windows is present and is addressed for the
bilinear model.

\subsection{GARCH Model}

We now consider the GARCH(1,1) model
\[
\begin{cases}
X_t=\sqrt{h_t}\,Z_t\\
h_t=\alpha_0+\alpha_1 X_{t-1}^2+\alpha_2 h_{t-1}
\end{cases}
\]
where $\bm{\alpha}=(.1,.8,.1)$ and
$Z_t\stackrel{\text{iid}}{\sim}\N(0,1)$.  The theoretical values of
the bispectrum are unknown, so they are approximated via simulation
over 500 realizations at a length of $10^5$ and averaging the four
estimators.

\begin{table}[H]
\centering
\begin{tabular}{|c|ccccc|}
  \hline
  $N=200$ & $|\hat{f}(0,0)|$ & Re$\hat{f}(2,1)$ & Im$\hat{f}(2,1)$ & $|\hat{f}(2,1)|$ & $T_6$ \\
  \hline
$\lambda_{\text{rpf}}$ & $\text{9.752e-04}$ & $\text{5.462e-05}$ & $\text{3.92e-05}$ & $\text{9.383e-05}$ & $\text{113.1}$ \\
$\lambda_{\text{rcf}}$ & $\text{1.038e-03}$ & $\text{5.800e-05}$ & $\text{4.391e-05}$ & $\text{1.019e-04}$ & $\text{115.1}$ \\
$\lambda_{\text{opt}} $ (f.p.) & $\text{6.580e-04}$ & $\text{4.345e-05}$ & $\text{3.182e-05}$ & $\text{7.527e-05}$ & $\text{110.1}$ \\
$\lambda_{\text{opt}} $ (s.p.) & $\text{3.849e-04}$ & $\text{3.488e-05}$ & $\text{5.112e-05}$ & $\text{8.600e-05}$ & $\text{125.1}$ \\
  \hline
  \hline
  $N=2000$ & $|\hat{f}(0,0)|$ & Re$\hat{f}(2,1)$ & Im$\hat{f}(2,1)$ & $|\hat{f}(2,1)|$ & $T_6$ \\
  \hline
$\lambda_{\text{rpf}}$ & $\text{2.411e-05}$ & $\text{2.916e-06}$ & $\text{1.555e-06}$ & $\text{4.471e-06}$ & $\text{7.317}$ \\
$\lambda_{\text{rcf}}$ & $\text{2.682e-05}$ & $\text{3.050e-06}$ & $\text{1.745e-06}$ & $\text{4.795e-06}$ & $\text{7.401}$ \\
$\lambda_{\text{opt}} $ (f.p.) & $\text{1.894e-05}$ & $\text{2.528e-06}$ & $\text{1.632e-06}$ & $\text{4.159e-06}$ & $\text{7.026}$ \\
$\lambda_{\text{opt}} $ (s.p.) & $\text{5.781e-05}$ & $\text{5.577e-06}$ & $\text{7.577e-06}$ & $\text{1.315e-05}$ & $\text{9.021}$ \\
  \hline
\end{tabular}
\caption{MSE estimates based on garch data for $N=200$ and
$N=2000$.}
\end{table}

For $N=200$, $\lambda_{\text{opt}}$ (s.p.) performed best at the
origin, but considerably worse in the composite criterion.  For the
larger $N$, the flat-top estimators and $\lambda_{\text{opt}}$
(f.p.) again performed significantly better than
$\lambda_{\text{opt}}$ (s.p.).

\subsection{Bilinear Model}

Finally, we consider the BL(1,0,1,1) bilinear model \cite{rao84}
\[
X_t=a X_{t-1}+bX_{t-1}Z_{t-1} + Z_t
\]
where $a=b=.4$ and $Z_t\stackrel{\text{iid}}{\sim}\N(0,1)$.  The
complete calculations of the bispectrum have been worked out in
\cite{rao84}, however the given equation for the bispectrum does not
match-up with the simulations.  Therefore theoretical values of the
bispectrum were computed through simulations as done in the GARCH
model.  The spectral density equation provided in \cite{rao84} is
correct and was used.

Whereas the previous three models had an optimal bandwidth of 1
throughout, the optimal bandwidths for the bilinear model is
typically much larger and depends on the evaluation criterion
considered.  The subscripted numbers represent the best and second
best bandwidth for each window (as deduced from simulation).

\begin{table}[H]
\centering
\begin{tabular}{|c|ccccc|}
  \hline
  $N=200$ & $|\hat{f}(0,0)|$ & Re$\hat{f}(2,1)$ & Im$\hat{f}(2,1)$ & $|\hat{f}(2,1)|$ & $T_6$ \\
  \hline
$\lambda_{\text{rpf}}$ & $\text{5.872}_{2,3}$ & $\text{5.421e-04}_{4,5}$ & $\text{1.008e-03}_{1,2}$ & $\text{1.55e-03}_{4,5}$ & $\text{806.8}_{2,5}$ \\
$\lambda_{\text{rcf}}$ & $\text{5.956}_{2,3}$ & $\text{6.005e-04}_{6,5}$ & $\text{1.073e-03}_{1,2}$ & $\text{1.673e-03}_{6,5}$ & $\text{817.0}_{7,6}$ \\
$\lambda_{\text{opt}} $ (f.p.) & $\text{4.401}_{2,1}$ & $\text{4.608e-04}_{4,3}$ & $\text{9.654e-04}_{1,2}$ & $\text{1.426e-03}_{4,3}$ & $\text{807.1}_{5,4}$ \\
$\lambda_{\text{opt}} $ (s.p.) & $\text{2.916}_{2,1}$ & $\text{3.926e-04}_{4,3}$ & $\text{8.623e-04}_{1,2}$ & $\text{1.255e-03}_{4,3}$ & $\text{791.4}_{5,4}$ \\
  \hline
    \hline
  $N=2000$ & $|\hat{f}(0,0)|$ & Re$\hat{f}(2,1)$ & Im$\hat{f}(2,1)$ & $|\hat{f}(2,1)|$ & $T_6$ \\
  \hline
$\lambda_{\text{rpf}}$ & $\text{1.755}_{4,3}$ & $\text{7.734e-05}_{4,6}$ & $\text{9.867e-05}_{1,2}$ & $\text{1.76e-04}_{4,6}$ & $\text{71.76}_{2,5}$ \\
$\lambda_{\text{rcf}}$ & $\text{1.891}_{4,3}$ & $\text{7.792e-05}_{6,7}$ & $\text{1.012e-04}_{1,2}$ & $\text{1.791e-04}_{6,7}$ & $\text{74.69}_{6,7}$ \\
$\lambda_{\text{opt}} $ (f.p.) & $\text{2.119}_{4,3}$ & $\text{6.282e-05}_{5,6}$ & $\text{9.443e-05}_{1,2}$ & $\text{1.572e-04}_{5,4}$ & $\text{71.01}_{6,7}$ \\
$\lambda_{\text{opt}} $ (s.p.) & $\text{1.322}_{4,3}$ & $\text{5.123e-05}_{5,6}$ & $\text{8.064e-05}_{1,2}$ & $\text{1.319e-04}_{5,4}$ & $\text{72.83}_{6,7}$ \\
 \hline
\end{tabular}
\caption{MSE estimates based on bilinear data for $N=200$ and
$N=2000$.} \label{tab:b}
\end{table}

For this model, $\lambda_{\text{opt}} $ (s.p.) performs better than
the other three, but with decreasing margins with increased $N$.
There is significant improvement of the flat-top estimators and
$\lambda_{\text{opt}}$ (f.p.) from $N=200$ to $N=2000$ making all
the estimators mostly equivalent.  The particularly good performance
of $\lambda_{\text{opt}}$ (s.p.) at the origin is due to a
fortuitous bandwidth selection under sensitive conditions; this is
addressed in more detail below.

There is somewhat of a discontinuity in optimal bandwidths for
$\lambda_{\text{rpf}}$ under the composite criterion as it jumps
from a best value of 2 to a second best value of 5.  A closer look
at the MSEs for each bandwidth from 1 to 8 further illustrates this.

\begin{table}[H]
\centering
\begin{tabular}{|c|cccccccc|}
  \hline
  $N=200$ & 1 &  2 &  3 &  4 & 5 & 6 & 7 & 8\\
  \hline
  $\lambda_{\text{rpf}}$ &  10200 &  331 &  851 &  442 &  423 &  457 &  440 & 454 \\
  $\lambda_{\text{rcf}}$ &  10200 &  471 &  985 &  561 &  458 &  457 &  454 & 464 \\
  $\lambda_{\text{opt}}$ &  6000 &   730 &  459 &  422 &  418 &  426 &  438 & 454 \\
  \hline
\end{tabular}
\caption{MSE estimates of $T_6$ with bandwidths one through ten and
$N=200$} \label{tab:b6}
\end{table}

We see that the bandwidth 2 is very good for the flat-top
lag-windows but very poor for $\lambda_{\text{opt}}$.  Moreover,
bandwidths 1 and 3 are extremely bad for the flat-top lag-window,
and any bandwidth larger than 3 is mostly equivalent among the
estimators.  In the bandwidth selection procedure only odd integer
bandwidths were selected since the last step of the procedure
generates the bandwidth from dividing an integer by $b=c=.51$. If
instead the parameter $c=.5$ is used, then only even integer
bandwidths would be produced by the algorithm.

The bispectrum corresponding to bilinear model resembles a hill
peaking at the origin \cite{rao84}.  This causes the choice of
bandwidth to be particularly delicate when estimating the origin.
The following table depicts this delicacy.

\begin{table}[H]
\centering
\begin{tabular}{|c|ccccccc|}
  \hline
  $N=200$ & 1 &  2 &  3 &  4 & 5 & 6 & 7 \\
  \hline
  $\lambda_{\text{rpf}}$    & 2.062 & {\bf 1.389} & 1.71 & 2.879 & 4.216 & 5.849 & 7.22   \\
  $\lambda_{\text{rcf}}$    & 2.062 & {\bf 1.390} & 1.864 & 3.207 & 4.848 & 6.502 & 8.078    \\
  $\lambda_{\text{opt}}$    & 1.823 & {\bf 1.445} & 2.013 & 3.13 & 4.448 & 5.733 & 6.852     \\
  \hline
\end{tabular}
\caption{MSE estimates at the origin with bandwidths one through
seven and $N=200$} \label{tab:b5}
\end{table}

We see that selecting any bandwidth besides 2, or possibly 3, leads
to a much larger mean square error.  The bispectrum, however, is
much flatter at points away from the origin, like the six interior
points used in the composite evaluation.  This causes the bandwidth
to be less sensitive to the choice of bandwidth when estimating an
interior value as seen in Table \ref{tab:b6} above.

The simulations up to this point mostly depict the strength of the
bandwidth selection procedure, and not the general asymptotic
optimality of the flat-top lag-window.  However, if we consider MSE
estimates for a fixed set of bandwidths, as in Table \ref{tab:b5},
the flat-top estimates perform better than $\lambda_{\text{opt}}$
which improves with $N$.  The following table demonstrates the
increased performance at $N=2000$.

\begin{table}[H]
\centering
\begin{tabular}{|c|cccccccc|}
  \hline
  $N=2000$ & 1 &  2 &  3 &  4 & 5 & 6 & 7 & 8\\
  \hline
  $\lambda_{\text{rpf}}$  & 2.029 & 0.9465 & 0.552 & {\bf 0.4687} & 0.6002 & 0.8262 & 1.029 & 1.237\\
  $\lambda_{\text{rcf}}$  & 2.029 & 0.9082 & 0.5074 & {\bf 0.4917} & 0.6821 & 0.9224 & 1.156 & 1.346   \\
  $\lambda_{\text{opt}}$  & 1.736 & 0.8919 & 0.5444 & {\bf 0.5267} & 0.6444 & 0.8001 & 0.9579 & 1.099   \\
  \hline
\end{tabular}
\caption{MSE estimates at the origin with bandwidths one through
seven and $N=200$} \label{tab:b7}
\end{table}

Further illustration of the optimality of the flat-top lag-windows
is provided in \cite{politis95} where second-order spectral density
estimation with flat-top lag-windows is addressed.

\subsection{Analysis of Bandwidth Procedures}

Histograms of the bandwidths produced by the procedures are provided
below.  A summary of their performance is tabulated in the following
table.

\ctable[caption={MSE of $\hat{M}/M-1$ for bandwidth selection
procedures (a)--(e)},botcap,label=resultsb,pos=H]{lccccccccccc}
{\tnote{Bandwidths 5 and 6 were selected as theoretical bandwidths
for procedure (a), but this is only approximate as the optimal
bandwidth varies.  True theoretical bandwidths can be inferred from
Table \ref{tab:b}. } } { \FL
 & \mc{2}{IID} & & \mc{2}{ARMA} & &\mc{2}{GARCH} & & \mc{2}{Bilinear\tmark}\NN
\cmidrule{2-3}\cmidrule{5-6}\cmidrule{8-9}\cmidrule{11-12} 
$N$  & 200  & 2000 &   & 200 & 2000 &    & 200  & 2000    &  & 200   & 2000   \ML
(a) & 3.18 & 0.792 &  & 1.54 & 0.248 &  & 6.36 & 0.968    &  & 0.413 & 0.182     \NN
(b) & 0.862 & 0.276 &  & 0.232 & .050 &  & 2.59 & 0.292   &  & 1.63 & 0.454     \NN
(c) & 2.71 & 0.900 &  & 0.866 & 0.142 &  & 4.05 & 0.552   &  & 0.633 & 0.362    \NN
(d) & 1.45 & 3.96 &  & 1.27 & 3.22 &  & 1.19 & 3.49       &  & 0.185 & 0.414       \NN
(e) & 4.66 & 12.0 &  & 4.04 & 9.36 &  & 4.22 & 9.82       &  & 0.0706 & 0.0394      \LL }

We see that the simple bandwidth selection algorithm is very
effective in producing accurate bandwidths that are consistent.  The
bandwidth selection procedure (a) can be seen to be quite accurate
from the histograms but tends to produce a few relatively large
bandwidths.  This error is compounded when squared error loss is
used to evaluate the performance.  The plug-in method with
second-order pilots on the other hand performs very poorly and does
not even appear consistent.

Histograms of the five bandwidth selection procedures are provided
in Appendix \ref{app:hist}. The histograms in the first three models
show a clear convergence of procedures (a) through (c) to the ideal
bandwidth 1, whereas the bandwidths from procedures (d) and (e) grow
with $N$.  The histograms for the bilinear model show a general
increase in $M$ with $N$ across each procedure.

\section{Conclusions}

Flat-top kernels in higher-order spectral density estimation is
shown to be asymptotically superior in terms of MSE to any other
finite-order kernel estimators.  In addition, a very simple
bandwidth selection algorithm is included that delivers ideal
bandwidths tailored to the flat-top estimators. If one chooses not
to adopt the infinite-order flat-top lag-window, then bandwidth
selection via the plug-in method with flat-top pilots demonstrates
greatly increased performance and should be used.  Finite-sample
simulations show these flat-top estimators were comparable with, and
in many cases outperforming, the popular second-order ``optimal''
lag-window estimator using the plug-in method with second-order
pilots for bandwidth selection.  Simulations show the estimation of
the bispectrum is quite sensitive to the choice of bandwidth, and
this paper delivers the first higher-order accurate bandwidth
selection procedures for the bispectrum.

\appendix

\section{Technical proofs}

\begin{lemma}
\label{lemma:expectation}
  The expectation of $\widehat{C}(\bm{\tau})$ is
  \[\E\left[\widehat{C}(\bm{\tau})\right]=\lrp{1-\frac{\gamma}{N}}C(\bm{\tau})+O\left(\frac{1}{N}\right).\]
\end{lemma}

\noindent{\sc Proof of Lemma \ref{lemma:expectation}.}  Let
$\bm{y}_t=\bm{x}_t-\bm{\mu}$, then
\begin{equation*}
\begin{split}
  \E\left[\widehat{C}(\bm{\tau})\right]&=\frac{1}{N}\sum_{t=1}^{N-\gamma}\E\left[
  \prod_{j=1}^{s}(x^{(a_j)}_{t-\alpha+\tau_j}-\bar{x}^{(a_j)})\right]\\
&=\frac{1}{N}\sum_{t=1}^{N-\gamma}\E\left[
  \prod_{j=1}^{s}\left((x^{(a_j)}_{t-\alpha+\tau_j}-\mu^{(a_j)})+(\mu^{(a_j)}-\bar{x}^{(a_j)})\right)\right]\\
&=\frac{1}{N}\sum_{t=1}^{N-\gamma}\E\left[
  \prod_{j=1}^{s}(y^{(a_j)}_{t-\alpha+\tau_j}-\bar{y}^{(a_j)})\right]\\
&=\frac{1}{N}\sum_{t=1}^{N-\gamma}C(\bm{\tau})+
  \sum_{\bm{\delta}\in\mathcal{V}}
  \E\left[\prod_{j=1}^{s}\left(y^{(a_j)}_{t-\alpha+\tau_j}\right)^{1-\delta_j}
  \left(\bar{y}^{(a_j)}\right)^{\delta_j}\right].\\
\end{split}
\end{equation*}
In the summation above, $\mathcal{V}$ denotes the set of all binary
$s$-tuples excluding the $s$-tuple $\{0,\ldots,0\}$; $\mathcal{V}$
has cardinality $2^s-1$. Let $\bm{\delta}\in\mathcal{V}$ and $\ell$
be its weight, i.e. $\ell=\sum_{j=1}^s\delta_j$.  Let us suppose,
w.l.o.g., that the first $\ell$ components of $\bm{\delta}$ are 1
and the rest 0. Then the term in the above summation corresponding
to this $\bm{\delta}$ can be written as
\begin{equation*}
\begin{split}
\E\left[\prod_{j=1}^{\ell}\bar{y}^{(a_j)}\prod_{j=\ell+1}^{s}y^{(a_j)}_{t-\alpha+\tau_j}
  \right]
&=\frac{1}{N^{\ell}}\sum_{u_1=1}^N\cdots\sum_{u_\ell=1}^N\E\left[
\prod_{i=1}^{\ell}y_{u_i}^{(a_i)}\prod_{j=\ell+1}^{s}y^{(a_j)}_{t-\alpha+\tau_j}\right]\\
&=\frac{1}{N^{\ell}}\sum_{u_1=1}^N\cdots\sum_{u_\ell=1}^N
C(u_1,\ldots,u_\ell,t-\alpha+\tau_{\ell+1},\ldots,t-\alpha+\tau_{s})\\
&=O\left(\frac{1}{N}\right).
\end{split}
\end{equation*}
The last equality follows from the absolute summability of $C(\bt)$.
Since for every $\bm{\delta}\in\mathcal{V}$ the expectation as above
is $O(N^{-1})$, the result follows.

\vspace{.5cm}

\noindent {\sc Proof of Theorem \ref{theorem:bias}.} Using Lemma
\ref{lemma:expectation} and property (iii) of the lag-window, the
expectation of $\hat{f}(\bm{\omega})$ can be expressed as
\[
\begin{split}
E[\hat{f}(\bm{\omega})]&=\frac{1}{(2\pi)^{s-1}}\sum_{\|\bt\|<N}
\left(\left(1-\frac{\gamma}{N}\right)C(\bt)+O\Big(\frac{1}{N}\Big)\right)
  \lambda_M(\bt)e^{-i\bt\cdot\bm{\omega}}\\
&=\frac{1}{(2\pi)^{s-1}}\sum_{\|\bt\|<N}\left(1-\frac{\gamma}{N}\right)C(\bt)
  \lambda_M(\bt)e^{-i\bt\cdot\bm{\omega}}+O\Big(\frac{M^{s-1}}{N}\Big).
\end{split}
\]

The bias of $\hat{f}(\bm{\omega})$ is
\begin{equation*}
\begin{split}
  E[\hat{f}(\bm{\omega})]-f(\bm{\omega})&=\frac{1}{(2\pi)^{s-1}}\sum_{\|\bt\|<N}
  \left(\lambda_M(\bt)\left(1-\frac{\gamma}{N}\right)-1\right)
  C(\bt)e^{-i\bt\cdot\bm{\omega}}\\
  &\qquad -\frac{1}{(2\pi)^{s-1}}\sum_{\|\bm{\tau}\|\ge N}
  C(\bt)e^{-i\bt\cdot\bm{\omega}}+O\Big(\frac{M^{s-1}}{N}\Big)\\
  &=\underbrace{\frac{1}{(2\pi)^{s-1}}\sum_{{\|\bm{\tau}\|<N}}
  \left(\lambda_M(\bt)-1\right)
  C(\bt)e^{-i\bt\cdot\bm{\omega}}}_{A_1}\\
  &\qquad \underbrace{-\frac{1}{(2\pi)^{s-1}N}\sum_{\|\bt\|<N}
  \gamma\,\lambda_M(\bt)\,
  C(\bt)e^{-i\bt\cdot\bm{\omega}}}_{A_2}\\
  &\qquad\underbrace{-\frac{1}{(2\pi)^{s-1}}\sum_{\|\bt\|\ge N}
  C(\bt)e^{-i\bt\cdot\bm{\omega}}}_{A_3}+O\Big(\frac{M^{s-1}}{N}\Big)
\end{split}
\end{equation*}
By the assumption on the summability of $C(\bt)$, $|A_3|$ can be
bounded as
\begin{equation*}
\begin{split}
  |A_3|\le \frac{1}{(2\pi)^{s-1}}\sum_{\|\bt\|\ge N}
  |C(\bt)|\le
  \frac{1}{(2\pi)^{s-1}N^{k}}\sum_{\|\bt\|\ge N}
  \|\bt\|^k|C(\bt)|=o\left(\frac{1}{N^k}\right)
\end{split}
\end{equation*}
Also,
\begin{equation*}
  |A_2|\le\frac{1}{(2\pi)^{s-1}N}\sum_{\|\bt\|<N}
  |\gamma|\,|C(\bt)|=O\left(\frac{1}{N}\right)
\end{equation*}
Now rewrite $A_1$ as
\begin{equation*}
\begin{split}
  A_1&=\underbrace{\frac{1}{(2\pi)^{s-1}}\sum_{\|\bt\|\le bM}
  \left(\lambda_M(\bt)-1\right)
  C(\bt)e^{-i\bt\cdot\bm{\omega}}}_{=0}\\
  &\qquad+\frac{1}{(2\pi)^{s-1}}\sum_{bM< \|\bt\|\le N}
  \left(\lambda_M(\bt)-1\right)
  C(\bt)e^{-i\bt\cdot\bm{\omega}}
\end{split}
\end{equation*}
\noindent {\sc Proof of ({\rm i}).}

Since $|\lambda(\bm{s})|\le 1$,
\begin{equation*}
\begin{split}
  |A_1|&\le \frac{2}{(2\pi)^{s-1}}\sum_{bM< \|\bt\|\le N}
  |C(\bt)|\le \frac{2}{(2\pi)^{s-1}(bM)^{k}}\sum_{bM< \|\bt\|\le N}
  \|\bt\|^k|C(\bt)|=o\left(\frac{1}{M^k}\right)
\end{split}
\end{equation*}
Equation (\ref{eq:bias}) now follows, and thus
$\MSE(\hat{f}(\bm{\omega}))\sim o(M^{-2k})+O(M^{s-1}/N)$.

\vspace{.5cm}

\noindent {\sc Proof of ({\rm ii}).}

  We have
  $\bias(\hat{f}(\bm{\omega}))=A_1+O\lrp{M^{s-1}/N}$, where under
  the assumptions of (ii),
\begin{equation*}
\begin{split}
  |A_1|\le \frac{2}{(2\pi)^{s-1}}\sum_{bM< \|\bt\|\le N}
  |C(\bt)|\le \frac{(2)\,D}{(2\pi)^{s-1}e^{dbM}}\sum_{bM< \|\bt\|\le N}
  e^{d(bM-\|\bt\|)}=O\left(e^{-dbM}\right).
\end{split}
\end{equation*}
Therefore $\MSE(\hat{f}(\bm{\omega}))\sim O(e^{-2dbM})+O(M^{s-1}/N)$
is asymptotically minimized when $M\sim A\log N$ where $A=1/(2db)$,
and (\ref{eq:mse2}) holds for all $A\ge 1/(2db)$.

\vspace{.2cm}

\noindent {\sc Proof of ({\rm iii}).}

We have
  $\bias(\hat{f}(\bm{\omega}))=A_1+O\lrp{M^{s-1}/N}$, but under
  the assumptions of (iii), $A_1=0$. Hence the bias and variance are
  $O(1/N)$.

\vspace{.5cm}

\noindent {\sc Proof of Theorem \ref{theorem:bandwidth}.} Let
$\bt_{\hat{m}}$ be any element of norm $\hat{m}$ for which
\begin{equation}
\label{eq:rhohat1} |\hat{\rho}(\bt_{\hat{m}})|>k\sqrt{\frac{\log
N}{N}}
\end{equation}
and let $\bt_{\hat{m}}'\in B_{\hat{m},\hat{m}+1}$, so that
$\hat{m}<\|\bt_{\hat{m}}'\|\le\hat{m}+1$, and
\begin{equation}
\label{eq:rhohat2} |\hat{\rho}(\bt_{\hat{m}}')|<k\sqrt{\frac{\log
N}{N}}
\end{equation}
 Equations (\ref{eq:alg}) and (\ref{eq:rho2}) give
\begin{equation}
\label{eq:rhohat4}
|\hat{\rho}(\bt_{\hat{m}})|=|\rho(\bt_{\hat{m}})|+o_p\left(\sqrt{\frac{\log\log
N}{N}}\right)
\end{equation}
In part (i), $\rho(\bt)\sim A\|\bt\|^{-d}$, so for any $\epsilon>0$,
we can find $\tau_0$ such that
\begin{equation}
\label{eq:rhobound1}
A(1-\epsilon)\|\bt\|^{-d}<\rho(\bt)<A(1+\epsilon)\|\bt\|^{-d}
\end{equation}
 when
$\|\bt\|>\tau_0$.  Similarly, for any $\epsilon>0$, there exists
$\tau_0$ large enough such that
\begin{equation}
\label{eq:mhat} (1-\epsilon)\hat{m}^{-d}<\,
\|\bt\|^{-d}\!<(1+\epsilon)\hat{m}^{-d}\qquad\text{for all }\bt\in
B_{\hat{m},\hat{m}+1}
\end{equation}
when $\hat{m}>\tau_0$.  Putting equations (\ref{eq:rhohat1}),
(\ref{eq:rhohat2}), (\ref{eq:rhohat4}), (\ref{eq:rhobound1}), and
(\ref{eq:mhat}) together gives, with high probability,
\begin{equation}
\label{eq:mhat2} A(1-\epsilon)^2\hat{m}^{-d}<c\sqrt{\frac{\log
N}{N}}<A(1+\epsilon)^2\hat{m}^{-d}
\end{equation}
up to $o_p(\sqrt{(\log\log N)/N})$, which is negligible as $N$ gets
large.  Equation (\ref{eq:mhat2}) is equivalent to
\[
\frac{\hat{m}}{(1+\epsilon)^2}<\frac{A^{1/ d}N^{1/2d}}{k^{1/d}(\log
N)^{1/2d}} <\frac{\hat{m}}{(1-\epsilon)^2}
\]
with high probability.  Therefore
\[
\hat{m}\stackrel{P}{\sim}\frac{A^{1/ d}N^{1/2d}}{k^{1/d}(\log
N)^{1/2d}}
\]
 The proof of
part (ii) is similar.

Now we prove part (iii).  Note that $\hat{m}>q$ only if
\begin{equation}
\label{eq:ma1}
 \max_{\bt\in
B_{q,q+a_N}}|\hat{\rho}(\bt)-\rho(\bt)|\ge k\sqrt{\frac{\log N}{N}}
\end{equation}
but since $C(\bt)=0$ when $\|\bt\|>q$, equation (\ref{eq:rho2}) then
shows
\begin{equation}
\label{eq:ma2} \max_{\bt\in
B_{q,q+a_N}}|\hat{\rho}(\bt)|=o_p\left(\sqrt{\frac{\log \log
N}{N}}\right)
\end{equation}
since $a_N=o(\log N)$.  The probability of (\ref{eq:ma1}) and
(\ref{eq:ma2}) happening simultaneously tends to zero, hence
$P(\hat{m}>q)\rightarrow 0$.  Now if $\hat{m}<q$ then
\[
B_{\hat{m},\hat{m}+a_N}|\hat{\rho}(\bt)|=|\rho(\bt)|+o_p\left(\sqrt{\frac{\log
\log N}{N}}\right)
\]
shows that (\ref{eq:alg}) must eventually be violated, hence
$P(\hat{m}<q)\rightarrow0$ and the result follows.

\noindent {\sc Proof of Theorem \ref{theorem:pilot}.}  Parts (ii)
and (iii) follow from Theorems \ref{theorem:bias} and
\ref{theorem:bandwidth} and the $\delta$-method; see
\cite{politis03} for more details.  For part (i), first note that
$\sum_{\bm{\tau}\in\bZ^{s-1}\backslash\{\bm{0}\}}\|\bm{\tau}\|^\alpha<\infty$
if and only if $\alpha>s-1$.  In order for
$\sum_{\bt\in\bZ^{s-1}}\|\bm{\tau}\|^{k+2}\,|C(\bm{\tau})|<\infty$,
for some $k\ge1$, $d$ must satisfy $d-k-2>s-1$ or $d>s+k+1\ge s+2$.
Now the results of Theorem 1 hold for
$\widehat{f_{\omega_i,\omega_j}}$ in replace of
$\hat{f}(\omega_1,\omega_2)$ for any positive integer $k<d-s-1$, in
particular for $k=\lceil d-s-2\rceil$.  From the proof of Theorem 1,
the bias is of order $o\left(1/M^k\right)$, and since the variance
is of smaller order, the result now follows from substituting $M$
with the rate $(N/\log N)^{1/2d}$ from Theorem 2 (i).

\section{Histograms}
\label{app:hist}

Below are histograms of the bandwidth selection procedures (a)
through (e) based on The top row in every Figure corresponds to
$N=200$ and the bottom row corresponds to $N=2000$.

\begin{figure}[H]
\centerline{\includegraphics{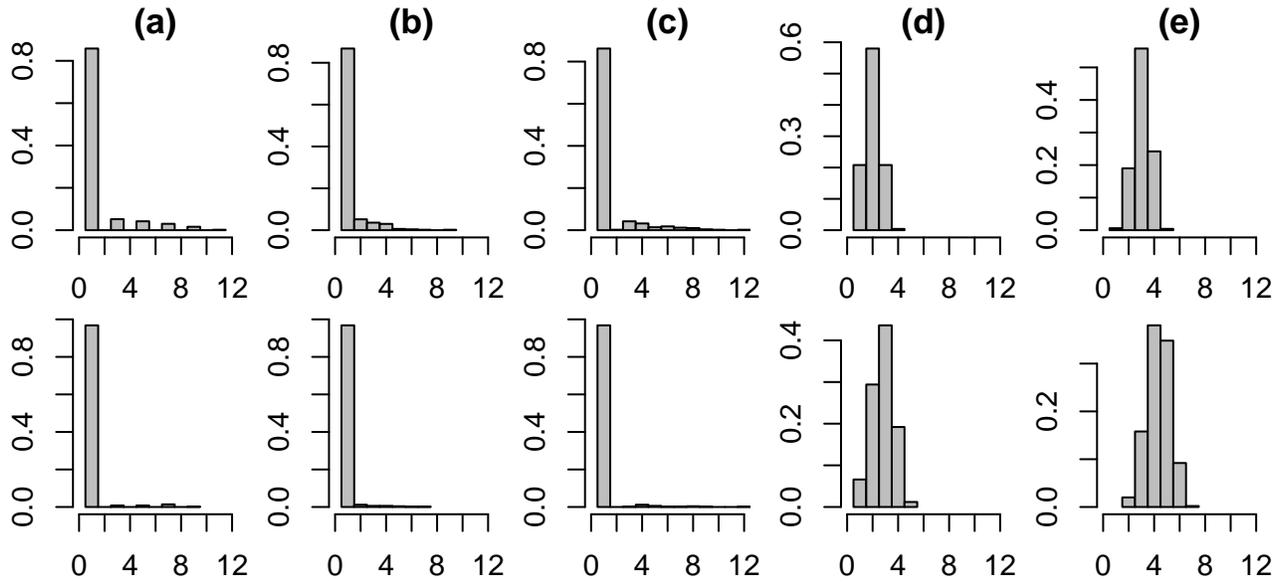}} \vspace{-.3cm}
\caption{Histograms based on iid data.} \label{fig:iid}
\end{figure}

\begin{figure}[H]
\centering
\includegraphics{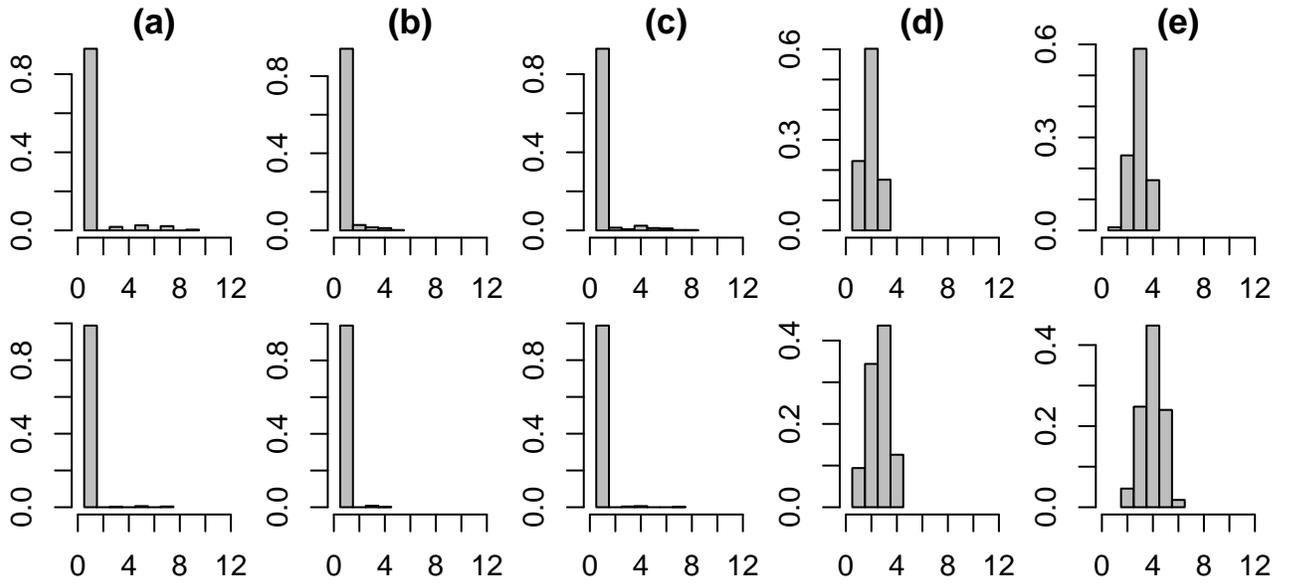} \vspace{-.3cm}
\caption{Histograms based on arma data.} \label{fig:arma}
\end{figure}

\begin{figure}[H]
\centerline{\includegraphics{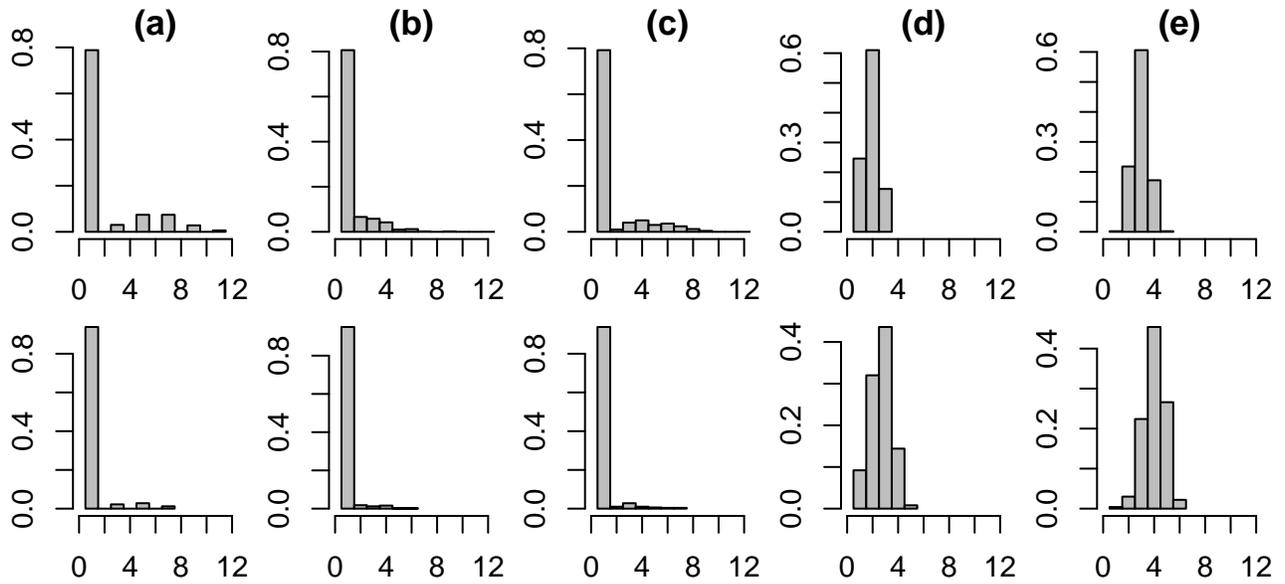}} \vspace{-.3cm}
\caption{Histograms based on garch data.}
\end{figure}

\begin{figure}[H]
\centering
\includegraphics{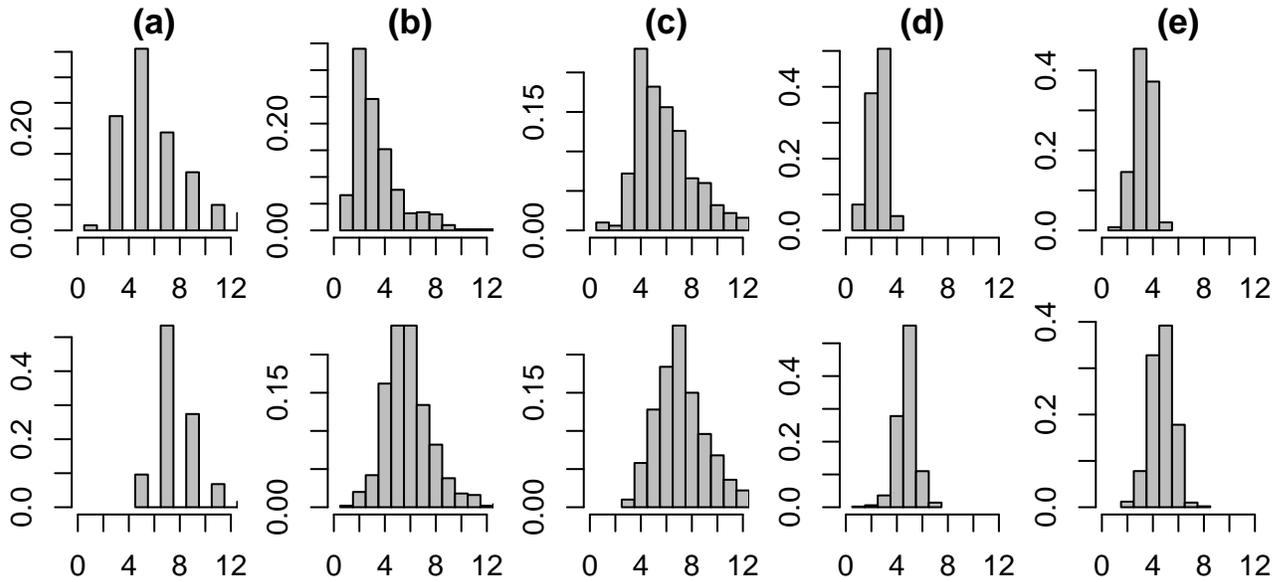} \vspace{-.3cm} \caption{Histograms based on
bilinear data.}
\end{figure}

\bibliography{flat_top}

\begin{thebibliography}{17}
\providecommand{\natexlab}[1]{#1}
\providecommand{\url}[1]{\texttt{#1}}
\expandafter\ifx\csname urlstyle\endcsname\relax
  \providecommand{\doi}[1]{doi: #1}\else
  \providecommand{\doi}{doi: \begingroup \urlstyle{rm}\Url}\fi

\bibitem[Brillinger and Rosenblatt(1967{\natexlab{a}})]{brillinger67a}
David~R. Brillinger and Murray Rosenblatt.
\newblock Asymptotic theory of estimates of {$k$}-th order spectra.
\newblock In \emph{Spectral Analysis Time Series (Proc. Advanced Sem., Madison,
  Wis., 1966)}, pages 153--188. John Wiley, New York, 1967{\natexlab{a}}.

\bibitem[Brillinger and Rosenblatt(1967{\natexlab{b}})]{brillinger67b}
David~R. Brillinger and Murray Rosenblatt.
\newblock Computation and interpretation of {$k$}-th order spectra.
\newblock In \emph{Spectral Analysis Time Series (Proc. Advanced Sem., Madison,
  Wis., 1966)}, pages 189--232. John Wiley, NEw York, 1967{\natexlab{b}}.

\bibitem[Brockmann et~al.(1993)Brockmann, Gasser, and Herrmann]{B93}
Michael Brockmann, Theo Gasser, and Eva Herrmann.
\newblock Locally adaptive bandwidth choice for kernel regression estimators.
\newblock \emph{J. Amer. Statist. Assoc.}, 88\penalty0 (424):\penalty0
  1302--1309, 1993.
\newblock ISSN 0162-1459.

\bibitem[B{\"u}hlmann(1996)]{B96}
Peter B{\"u}hlmann.
\newblock Locally adaptive lag-window spectral estimation.
\newblock \emph{J. Time Ser. Anal.}, 17\penalty0 (3):\penalty0 247--270, 1996.
\newblock ISSN 0143-9782.

\bibitem[Hall and Marron(1988)]{marron88}
Peter Hall and J.~S. Marron.
\newblock Choice of kernel order in density estimation.
\newblock \emph{Ann. Statist.}, 16\penalty0 (1):\penalty0 161--173, 1988.
\newblock ISSN 0090-5364.

\bibitem[Hinich(1982)]{Hinich82}
Melvin~J. Hinich.
\newblock Testing for {G}aussianity and linearity of a stationary time series.
\newblock \emph{J. Time Ser. Anal.}, 3\penalty0 (3):\penalty0 169--176, 1982.
\newblock ISSN 0143-9782.

\bibitem[Jammalamadak et~al.(2006)Jammalamadak, Rao, and Terdik]{J06}
S.R. Jammalamadak, T.S. Rao, and Gy{\"o}rgy Terdik.
\newblock Higher order cumulants of random vectors and applications to
  statistical inference and time series.
\newblock \emph{Sankhy\=a}, 68\penalty0 (2):\penalty0 326--356, 2006.

\bibitem[Jones et~al.(1996)Jones, Marron, and Sheather]{marron96}
M.~C. Jones, J.~S. Marron, and S.~J. Sheather.
\newblock A brief survey of bandwidth selection for density estimation.
\newblock \emph{J. Amer. Statist. Assoc.}, 91\penalty0 (433):\penalty0
  401--407, 1996.
\newblock ISSN 0162-1459.

\bibitem[Leadbetter et~al.(1983)Leadbetter, Lindgren, and
  Rootz{\'e}n]{Leadbetter83}
M.~R. Leadbetter, Georg Lindgren, and Holger Rootz{\'e}n.
\newblock \emph{Extremes and related properties of random sequences and
  processes}.
\newblock Springer Series in Statistics. Springer-Verlag, New York, 1983.
\newblock ISBN 0-387-90731-9.

\bibitem[Lii and Rosenblatt(1990{\natexlab{a}})]{lii90cov}
K.~S. Lii and M.~Rosenblatt.
\newblock Cumulant spectral estimates: bias and covariance.
\newblock In \emph{Limit theorems in probability and statistics (P\'ecs,
  1989)}, volume~57 of \emph{Colloq. Math. Soc. J\'anos Bolyai}, pages
  365--405. North-Holland, Amsterdam, 1990{\natexlab{a}}.

\bibitem[Lii and Rosenblatt(1990{\natexlab{b}})]{lii90norm}
K.~S. Lii and M.~Rosenblatt.
\newblock Asymptotic normality of cumulant spectral estimates.
\newblock \emph{J. Theoret. Probab.}, 3\penalty0 (2):\penalty0 367--385,
  1990{\natexlab{b}}.
\newblock ISSN 0894-9840.

\bibitem[Politis(2003)]{politis03}
Dimitris~N. Politis.
\newblock Adaptive bandwidth choice.
\newblock \emph{J. Nonparametr. Stat.}, 15\penalty0 (4-5):\penalty0 517--533,
  2003.
\newblock ISSN 1048-5252.

\bibitem[Politis and Romano(1995)]{politis95}
Dimitris~N. Politis and Joseph~P. Romano.
\newblock Bias-corrected nonparametric spectral estimation.
\newblock \emph{J. Time Ser. Anal.}, 16\penalty0 (1):\penalty0 67--103, 1995.
\newblock ISSN 0143-9782.

\bibitem[Saito and Tanaka(1985)]{Saito85}
Keiichi Saito and Tomoharu Tanaka.
\newblock Exact analytic expression for gabr-rao's optimal bispectral
  two-dimensional lag window.
\newblock \emph{J. Nucl. Sci. Technol.}, 22\penalty0 (12):\penalty0 1033--1035,
  1985.

\bibitem[Subba~Rao and Gabr(1980)]{gabr80}
T.~Subba~Rao and M.~M. Gabr.
\newblock A test for linearity of stationary time series.
\newblock \emph{J. Time Ser. Anal.}, 1\penalty0 (2):\penalty0 145--158, 1980.
\newblock ISSN 0143-9782.

\bibitem[Subba~Rao and Gabr(1984)]{rao84}
T.~Subba~Rao and M.~M. Gabr.
\newblock \emph{An introduction to bispectral analysis and bilinear time series
  models}, volume~24 of \emph{Lecture Notes in Statistics}.
\newblock Springer-Verlag, New York, 1984.
\newblock ISBN 0-387-96039-2.

\bibitem[Subba~Rao and Terdik(2003)]{rao03}
T.~Subba~Rao and Gy. Terdik.
\newblock On the theory of discrete and continuous bilinear time series models.
\newblock In \emph{Stochastic processes: modelling and simulation}, volume~21
  of \emph{Handbook of Statist.}, pages 827--870. North-Holland, Amsterdam,
  2003.

\end{thebibliography}

\end{document}